\documentstyle[amscd,amssymb,verbatim,12pt]{amsart}
\newcommand{\tr}{\operatorname{tr}}

\newcommand{\coh}{\operatorname{coh}}

\renewcommand{\mod}{\operatorname{mod}}

\newcommand{\Tor}{\operatorname{Tor}}

\newcommand{\OO}{{\cal O}}

\newcommand{\DD}{{\cal D}}

\newcommand{\BB}{{\cal B}}

\newcommand{\SL}{\operatorname{SL}}

\newcommand{\lan}{\langle}
\newcommand{\ran}{\rangle}
\newcommand{\Coh}{\operatorname{Coh}}

\newcommand{\CC}{{\cal C}}
\newcommand{\TT}{{\cal T}}
\newcommand{\MM}{{\cal M}}

\renewcommand{\P}{{\Bbb P}}

\newcommand{\Pic}{\operatorname{Pic}}

\newcommand{\de}{\delta}
\newcommand{\eps}{\epsilon}

\newcommand{\D}{{\cal D}}

\numberwithin{equation}{section}

\newtheorem{thm}{Theorem}[section]
\newtheorem{prop}[thm]{Proposition}
\newtheorem{prop-def}[thm]{Proposition-Definition}
\newtheorem{lem}[thm]{Lemma}
\newtheorem{cor}[thm]{Corollary}
\newenvironment{rem}{\vspace{3mm}\noindent
{\bf Remark.}}{\vspace{3mm}}
\newenvironment{defi}{\vspace{3mm}\noindent
{\bf Definition.}}{\vspace{3mm}}

\newcommand{\Pf}{\noindent {\it Proof}}
\newcommand{\id}{\operatorname{id}}

\newcommand{\ov}{\overline}

\newcommand{\rk}{\operatorname{rk}}

\renewcommand{\AA}{{\cal A}}
\newcommand{\FF}{{\cal F}}

\newcommand{\SS}{{\cal S}}

\newcommand{\Hom}{\operatorname{Hom}}
\newcommand{\Ext}{\operatorname{Ext}}

\renewcommand{\a}{\alpha}

\newcommand{\om}{\omega}

\newcommand{\la}{\lambda}
\newcommand{\th}{\theta}
\newcommand{\C}{{\Bbb C}}
\newcommand{\R}{{\Bbb R}}
\newcommand{\Z}{{\Bbb Z}}

\newcommand{\La}{\Lambda}

\newcommand{\wt}{\widetilde}
\newcommand{\ot}{\otimes}

\newcommand{\sub}{\subset}
\newcommand{\ed}{\qed\vspace{3mm}}

\newcommand{\Hol}{\operatorname{Hol}}

\newcommand{\vareps}{\varepsilon}
\newcommand{\Aut}{\operatorname{Aut}}

\title{Holomorphic bundles on 2-dimensional noncommutative toric orbifolds}
\author{A. Polishchuk}
\thanks{Supported in part by NSF grant}

\begin{document}
\begin{abstract} We define the notion of a holomorphic bundle on the noncommutative toric
orbifold $T_{\th}/G$
associated with an action of a finite cyclic group $G$ on an irrational rotation algebra.
We prove that the category of such holomorphic bundles is abelian and its derived category
is equivalent to the derived category of modules over a finite-dimensional algebra $\La$.
As an application we finish the computation of $K_0$-groups of the crossed 
product algebras describing the above orbifolds initiated in \cite{K}, \cite{W-Fourier}, \cite{W-Ktheory}, \cite{BW1} and \cite{BW2}.
Also, we describe a torsion pair in the category of $\La$-modules, such that the tilting with
respect to this torsion pair gives the category of holomorphic bundles on $T_{\th}/G$.
\end{abstract}
\maketitle

\bigskip

\centerline{\sc Introduction}

\medskip

Let $A_{\th}$ be the algebra of smooth functions on the noncommutative $2$-torus 
$T_{\th}$ associated with an irrational real number $\th$. Recall that its elements are expressions of
the form $\sum_{m,n}a_{m,n}U_1^mU_2^n$, where the coefficients $(a_{m,n})_{(m,n)\in\Z^2}$ 
rapidly decrease at infinity, and $U_1$ and $U_2$ satisfy the relation
$$U_1U_2=\exp(2\pi i\th)U_2U_1.$$
It is convenient to denote 
$U_v=\exp(-\pi i mn\th)U_1^mU_2^n$ 
for $v=(m,n)\in\Z^2$, so that 
$$U_vU_w=\exp(\pi i\th\det(v,w))U_{v+w}.$$
There is a natural action of $\SL_2(\Z)$ on $A_{\th}$ such that
the matrix $g$ acts by the automorphism $U_v\mapsto U_{gv}$. 
Hence, for a finite subgroup $G\sub\SL_2(\Z)$ we can
consider the crossed product algebra $B_{\th}=A_{\th}*G$.

The simplest case is when $G=\Z/2\Z$ generated by $-\id\sub\SL_2(\Z)$
acting on $A_{\th}$ by the so called {\it flip automorphism} .
In this case the algebra $B_{\th}$ was studied in the papers \cite{BEEK1}, \cite{BK}, \cite{K} and \cite{W}.
In particular, it is known that it is simple, has a unique tracial state, and is an AF-algebra.
Also its $K$-theory has been computed: one has $K_0(B_{\th})=\Z^6$ and $K_1(B_{\th})=0$.
However, there are three more examples of finite subgroups $G\sub\SL_2(\Z)$ for
which the situation is not so well understood. Namely, we can take $G=\Z/3\Z$ generated by
$\left(\matrix -1 & 1\\ -1 & 0 \endmatrix\right)$; or
$G=\Z/4\Z$ generated by the ``Fourier" matrix $\left(\matrix 0 & 1 \\ -1 & 0\endmatrix\right)$;
or $G=\Z/6\Z$ generated by $\left(\matrix -1 & 1 \\ 1 & 0 \endmatrix\right)$.
In this paper we compute $K_0(B_{\th})$ in all these cases using holomorphic vector bundles
on the corresponding noncommutative orbifolds.
 
By a {\it vector bundle} on the noncommutative toric orbifold
$T_{\th}/G$ we mean a finitely generated projective right $B_{\th}$-module.
We want to define what is a holomorphic structure on such a vector bundle.
As in \cite{PS}, \cite{P},
let us consider a complex structure on $T_{\th}$ associated with a complex number 
$\tau\in\C\setminus\R$. It is given by a derivation
$$\de:A_{\th}\to A_{\th}:\sum_{m,n}a_{m,n}U_1^mU_2^n\mapsto 2\pi i\sum_{m,n}
(m\tau+n)U_1^mU_2^n$$
of $A_{\th}$ that we view as an analogue of the $\ov{\partial}$-operator.
To descend this structure to the orbifold $T_{\th}/G$ we have to impose some compatibility
between the action of $G$ and $\de$. More precisely, we assume that there exists
a character $\vareps:G\to\C^*$ such that the following relation holds:
\begin{equation}\label{gdelta}
g\de=\vareps(g)\de g
\end{equation}
for all $g\in G$. 
For example, for $G=\Z/2\Z$ acting by the flip automorphism $\vareps$ is the unique nontrivial
character of $\Z/2\Z$.
In three other cases such a relation exists for a special choice of $\tau$. Namely,
let us identify $G=\Z/m\Z$ with the subgroup of $m$-th roots of unity in $\C^*$ and let
$G$ act on $\C$ by multiplication. Then we can choose $\tau$ in such a way that the lattice $\Z\tau+\Z$
is $G$-invariant: for $m=4$ we take $\tau=i$, while for $m=3$ and $m=6$ we 
take $\tau=(1+i\sqrt{3})/2$. Note that the embedding of $G$ into $\SL_2(\Z)$ is induced by its
action on the basis $(\tau, 1)$ of $\Z\tau+\Z$. 
Then \eqref{gdelta} will hold with $\vareps(g)=g^{-1}\in\C^*$.

Recall that in \cite{PS}, \cite{P} we studied the category $\Hol(T_{\th,\tau})$ of 
holomorphic bundles on $T_{\th}$. By definition, these are pairs $(P,\ov{\nabla})$ consisting of a finitely
generated projective right $A_{\th}$-module $P$ and an operator
$\ov{\nabla}:P\to P$ satisfying the Leibnitz identity 
$$\ov{\nabla}(f\cdot a)=f\cdot\de(a)+\ov{\nabla}(f)\cdot a,$$
where $f\in P$, $a\in A_{\th}$.
Now we extend $\de$ to a {\it twisted derivation} $\wt{\de}$ of $B_{\th}$ by setting
$$\wt{\de}(\sum_{g\in G} a_g g)=\sum_{g\in G} \eps(g)\de(a_g)g,$$
where $a_g\in A_{\th}$ for $g\in G$.
This extended map satisfies the twisted Leibnitz identity
$$\wt{\de}(b_1b_2)=b_1\wt{\de}(b_2)+\wt{\de}(b_1)\kappa(b_2),$$
where $\kappa$ is the automorphism of $B_{\th}$ given by
$\kappa(\sum_{g\in G} a_g g)=\sum_{g\in G} \eps(g) a_g g.$
We define a {\it holomorphic structure} on a vector bundle $P$ on $T_{\th}/G$ as
an operator $\ov{\nabla}:P\to P$ satisfying the similar twisted Leibnitz identity
$$\ov{\nabla}(f\cdot b)=f\cdot\wt{\de}(b)+\ov{\nabla}(f)\cdot \kappa(b),$$
where $f\in P$, $b\in B_{\th}$. By definition, a {\it holomorphic bundle} is a pair $(P,\ov{\nabla})$ consisting of a vector bundle $P$ equipped with a holomorphic structure $\ov{\nabla}$. 
One can define morphisms between holomorphic bundles in a natural way, so we obtain the category $\Hol(T_{\th,\tau}/G)$ of holomorphic bundles.

Recall that the combined results of \cite{PS} and \cite{P} imply that
the category $\Hol(T_{\th,\tau})$ is abelian and one has an equivalence of bounded derived categories
$$D^b(\Hol(T_{\th,\tau}))\simeq D^b(\Coh(E)),$$
where $\Coh(E)$ is the category of coherent sheaves on the elliptic curve $E=\C/(\Z+\Z\tau)$.
Furthermore, the image of the abelian category
$\Hol(T_{\th,\tau})$ in the derived category $D^b(\Coh(E))$ can be described as the
heart of the tilted $t$-structure associated with a certain torsion pair in $\Coh(E)$ (depending
on $\th$).
Our main result is a similar explicit description of the category of holomorphic bundles on 
$T_{\th,\tau}/G$, where $G=\Z/m\Z\sub\SL_2(\Z)$ with $m\in\{2,3,4,6\}$.

\begin{thm}\label{mainthm}
The category $\Hol(T_{\th,\tau}/G)$ is abelian and one has an equivalence of
bounded derived categories
$$D^b(\Hol(T_{\th,\tau}/G))\simeq D^b(\mod-\La),$$
where $\mod-\La$ is the category of finite-dimensional right modules over the 
algebra $\La=\C Q/(I)$ of paths in a quiver $Q$ without cycles
with quadratic relations $I$. The number of vertices of $Q$ is equal to
$6, 8, 9$ or $10$ for $m=2,3,4$ or $6$, respectively. 
\end{thm}

The precise description of the algebra $\La$ will be given in section \ref{semiort-sec}.
It is derived equivalent to one of  canonical tubular algebras considered by Ringel in \cite{R}.
Furthermore, the image of $\Hol(T_{\th,\tau}/G)$ in $D^b(\mod-\La)$ corresponds
to the tilted $t$-structure for a certain explicit torsion pair in $\mod-\La$ depending on $\th$
that we will describe in Theorem \ref{tilting-thm}. 

We prove Theorem \ref{mainthm} in two steps: first, we relate holomorphic bundles on
$T_{\th,\tau}/G$ to the derived category of $G$-equivariant sheaves on
the elliptic curve $E=\C/(\Z\tau+\Z)$, and then we construct a derived equivalence with
right modules over the algebra $\La$. The second step is actually well known and works for arbitrary weighted projective curves considered in \cite{GL}. We present an alternative
derivation working directly with equivariant sheaves. It is based on the semiorthogonal
decomposition of the category of $G$-equivariant sheaves associated with a ramified
$G$-covering of smooth curves (see Theorem \ref{semiort-thm}).

Combining Theorem \ref{mainthm} with the results of \cite{W-Fourier} and \cite{BW1}
we derive the following result.

\begin{thm}\label{K-thm} 
One has $K_0(B_{\th})\simeq\Z^r$, where $r=6,8,9$ or $10$ for $G=\Z/m\Z$ with $m=2,3,4$ or $6$,
respectively.
\end{thm}

Note that for $G=\Z/2\Z$ this was known (see \cite{K}). 
For $G=\Z/4\Z$ and $G=\Z/6\Z$ this was proved for $\th$ in a dense $G_{\delta}$-set  
(see \cite{W-Ktheory} and \cite{BW2}). The case of $G=\Z/4\Z$ and arbitrary irrational $\th$
was done by Lueck, Phillips and Walters in 2003 (unpublished).
Our proof shows in addition that the natural forgetful map
$$K_0(\Hol(T_{\th,\tau}/G))\to K_0(B_{\th})$$ 
is, in fact, an isomorphism and identifies the
positive cones in these groups. 

\noindent
{\it Acknowledgment}. I am grateful to Julian Buck, Igor Burban, Pavel Etingof, Helmut Lenzing, Tony Pantev, Chris Phillips, Olivier Schiffman and Samuel Walters for useful discussions.
During the first stages of work on this paper the author enjoyed hospitality of the Mathematisches Forschungsinstitut Oberwolfach. I'd like to thank this institution 
for providing excellent working conditions and a stimulating atmosphere.  

\section{Derived categories of $G$-sheaves}
\label{quasi-sec}

\subsection{Generalities on $G$-sheaves}

Let $G$ be a finite group acting on an algebraic variety $X$ over a field $k$ of characteristic zero. Then we can consider the category
$\Coh_G(X)$ of $G$-equivariant coherent sheaves. We denote its bounded derived category
by $D^b_G(X)$. It is equivalent to the full subcategory in the bounded derived category
of $G$-equivariant quasicoherent sheaves consisting of complexes with coherent cohomology
(see Corollary 1 in \cite{Bezr}). We refer to \cite{BKR}, Section 4, for a more detailed discussion of
this category and restrict ourself to several observations. Below we will use the term
{\it $G$-sheaf} to denote a $G$-equivariant coherent sheaf.

There is a natural forgetting functor from $D^b_G(X)$ to $D^b(X)$, the usual derived category
of coherent sheaves on $X$. For equivariant complexes of sheaves $F_1$ and $F_2$ we denote 
by $\Hom_G(F_1,F_2)$ (resp., $\Hom(F_1,F_2)$) morphisms between these objects in the former (resp., latter) category. There is a natural action of $G$ on $\Hom(F_1,F_2)$ and one has
$$\Hom_G(F_1,F_2)\simeq\Hom(F_1,F_2)^G.$$
In particular, the cohomological dimension of $\Coh_G(X)$ is at most that of $\Coh(X)$.
Let us also set $\Hom^i_G(F_1,F_2)=\Hom_G(F_1,F_2[i]$. If $X$ is a smooth projective variety
over $k$ then we can define the bilinear form
$\chi_G(\cdot,\cdot)$ on $K_0(D^b_G(X))$ by setting
$$\chi_G(F_1,F_2)=\sum_{i\in\Z}(-1)^i\dim\Hom_G^i(F_1,F_2).$$

Many natural constructions with sheaves carry easily to $G$-equivariant setting.
For example, the tensor product of $G$-sheaves is defined.
Also if $\rho$ is a representation of $G$ then there is a natural tensor product operation
$F\mapsto F\otimes\rho$ on $G$-sheaves.
If $Y$ is another variety equipped with the action of $G$ and if
$f:X\to Y$ is a $G$-equivariant morphism then there are natural functors of
push-forward and pull-back:
$$f_*:\Coh_G(X)\to\Coh_G(Y),\  f^*:\Coh_G(Y)\to\Coh_G(X),$$ 
We can also consider the derived functor $Rf_*:D^b_G(X)\to D^b_G(Y)$ and 
if $Y$ is smooth or $f$ is flat, the derived functor
$Lf^*:D^b_G(Y)\to D^b_G(X)$ (when $f$ is flat we denote it simply by $f^*$).
The pair $(Lf^*,Rf_*)$ satisfies the usual adjunction property.

If $X$ is smooth and projective then the category $D^b_G(X)$ is also equipped with
the Serre duality of the form
$$\Hom_G(F_1,F_2)^*\simeq\Hom_G(F_2,F_1\otimes\om_X[\dim X]),$$
where the canonical bundle $\om_X$ is equipped with the natural $G$-action.

The following observation will be useful to us.

\begin{lem}\label{eq-sh-lem} 
Let $X$ be a smooth curve equipped with an action of a finite group $G$.
Then the category $D^b_G(X)$ is equivalent to the category of $G$-equivariant
objects in $D^b(X)$, i.e., the category of data $(F,\phi_g)$, where $F\in D^b(X)$ and
$\phi_g:g^*F\wt{\to} F$, $g\in G$, is a collection of isomorphisms satisfying the natural compatibility
condition.
\end{lem}

\Pf . Note that there is a natural functor from $D^b_G(X)$ to the category of $G$-objects
in $D^b(X)$. It is easy to see that it is fully faithful, so the only issue is to check that it is essentially surjective. The proof is based on the fact that 
every object $F\in D^b(X)$ is isomorphic to the direct sum of its cohomology sheaves:
$F\simeq \oplus_n H^nF[-n]$.
A $G$-structure on $F$ is given by a compatible collection of isomorphisms $\phi=(\phi_g)$, where
$\phi_g:\oplus_n g^*F\to F$. Note that the only nontrivial components of $\phi_g$
with respect to the above direct sum decompositions are maps $g^*H^nF[-n]\to H^nF[-n]$
and  $g^*H^nF[-n]\to H^{n-1}F[-n+1]$. Let $\phi^0=(\phi^0_g)$ be the $G$-structure
on $F$ given by the components of $\phi$ of the first kind (i.e., by the diagonal components).
Since the decomposition of $F$ into the direct sum of cohomology sheaves is compatible with
$\phi^0$, it suffices to find an isomorphism of $G$-objects  
\begin{equation}\label{F-isom}
(F,\phi^0)\simeq (F,\phi).
\end{equation}
Let us write $\phi_g=a_g\circ\phi^0_g$, where $a_g\in\Aut(F)$. Note that
$a_g$ belongs to the abelian subgroup
$$A:=\oplus_n \Hom(H^nF[-n],H^{n-1}F[-n+1])\sub\Aut(F)$$
of "upper-triangular" automorphisms with identities as diagonal entries.
It is easy to check that the compatibility condition on the data $\phi$ amounts to
the $1$-cocycle equation for $a_g$, where $G$ acts on $A$ in a natural way.
On the other hand, existence of an isomorphism \eqref{F-isom} is equivalent
to $g\mapsto a_g$ being a coboundary. It remains to note that
$H^1(G,A)=0$ since $A$ is a vector space over a field of characteristic zero.
\ed

%\begin{rem} Can be generalized to an arbitrary $G$-action on an abelian category of cohomological
%dimension $\le 1$.
%\end{rem}

\subsection{Semiorthogonal decomposition associated with a Galois covering}
\label{semiort-sec}

Recall (see \cite{BKap}, \cite{BO}) that a {\it semiorthogonal decomposition} of a triangulated category $\AA$
is an ordered collection $(\BB_1,\ldots,\BB_r)$ of full triangulated subcategories in $\AA$ such that 
$\Hom(B_i,B_j)=0$ whenever $B_i\in\BB_i$ and $B_j\in\BB_j$, where $i>j$, and the
subcategories $\BB_1,\ldots,\BB_r$ generate $\AA$.
In this case we write
$$\AA=\lan \BB_1,\ldots,\BB_r\ran.$$

Semiorthogonal decompositions are related to admissible triangulated subcategories.
For a subcategory $\BB\sub\AA$ let us denote by $\BB^{\perp}$ the {\it right orthogonal} of $\BB$,
i.e., the full subcategory of $\AA$ consisting of all $C$ such that $\Hom(B,C)=0$ for
all $B\in\BB$.
A triangulated subcategory $\BB\sub\AA$ is called {\it right admissible} if
for every $X\in\AA$ there exists a distinguished triangle
$B\to X\to C\to\ldots$ with $B\in\BB$ and $C\in\BB^{\perp}$.
Thus, a right admissible subcategory $\BB\sub\AA$ gives rise to a semiorthogonal decomposition
$$\AA=\lan\BB^{\perp},\BB\ran.$$
Similarly, one can define the left orthogonal and left admissibility of a subcategory.

We are going to use also some results from the theory of exceptional collections (see 
\cite{Bondal}, \cite{R}). Let $k$ be a field.
Recall that an object $E$ in a $k$-linear triangulated subcategory $\AA$ is {\it exceptional} if
$\Hom^i(E,E)=0$ for $i\neq 0$ and $\Hom(E,E)=k$.
An {\it exceptional collection} in $\AA$ is a collection of exceptional objects
$(E_1,\ldots,E_n)$ such that $\Hom^*(E_i,E_j)=0$ for $i>j$. 
A triangulated subcategory $\lan E_1,\ldots,E_n\ran$
generated by an exceptional collection is known to be left and right admissible 
(see \cite{Bondal}, Theorem 3.2). In the case when $\lan E_1,\ldots,E_n\ran=\AA$ we will say
that $(E_1,\ldots,E_n)$ is a {\it full exceptional collection}.
An exceptional collection $(E_1,\ldots,E_n)$ is {\it strong} if $\Hom^a(E_i,E_j)=0$ for $a\neq 0$
and all $i,j$. If an exceptional collection is full and strong then 
$E=\oplus_{i=1}^n E_i$ is a tilting object in $\AA$, i.e., the functor
$X\mapsto R\Hom(E,X)$ gives an equivalence between $\AA$ and $D^b(\mod-A)$
(provided $\AA$ satisfies some natural finiteness assumptions and is framed, see
\cite{BKap2}).

Let $\pi:X\to Y$ be a ramified Galois covering with Galois group $G$, where $X$ and $Y$ are 
smooth projective curves over an algebraically closed field $k$ of characteristic zero. 
In other words, a finite group $G$ acts effectively on $X$ and $Y=X/G$.
We are going to construct a semiorthogonal decomposition of the derived category
of $G$-sheaves $D^b_G(X)$ with $D^b(Y)$ as one of the pieces.
Let $D_1,\ldots,D_r$ be all special fibers of $\pi$ equipped with the reduced scheme structure
and let $m_1,\ldots,m_r$ be the corresponding multiplicities. 
Let us also fix a point $p_i\in D_i$ for each $i=1,\ldots,r$, and let $G_i\sub G$ be the stabilizer subgroup of $p_i$. Then we have $G$-equivariant isomorphisms $D_i\simeq G/G_i$.
Hence, the category of $G$-sheaves on $D_i$ is equivalent to finite dimensional representations
of $G_i$. For every character $\zeta$ of $G_i$ we denote by $\zeta_{D_i}$ the corresponding
$G$-sheaf on $D_i$.
Note that $G_i$ is a cyclic group of order $m_i$. Moreover, the representation of $G_i$ on
$\om_X|_{p_i}$ allows to identify $G_i$ with the group of $m_i$-th roots of unity in
such a way that it acts on $\om_X|_{p_i}$ via the standard character. Thus, we have
an isomorphism of $G$-sheaves
$$\om_X|_{D_i}\simeq \zeta(i)_{D_i},$$
where $\zeta(i)$ is a generator of the character group $\hat{G}_i$.

\begin{thm}\label{semiort-thm} 
(i) The natural functor $\pi^*:D^b(Y)\to D^b_G(X)$ is fully faithful.

\noindent
(ii) For every $i=1,\ldots,r$, the collection of $G$-sheaves on $X$
\begin{equation}\label{exc-col1}
(\OO_{(m_i-1)D_i},\ldots,\OO_{2D_i},\OO_{D_i})
\end{equation}
is exceptional. Let $\BB_i$ be the triangulated subcategory in $D^b_G(X)$ generated
by this exceptional collection. Then $\BB_i$ and $\BB_j$ are mutually orthogonal
for $i\neq j$, i.e., $\Hom(B_i,B_j)=\Hom(B_j,B_i)=0$ for all $B_i\in\BB_i$ and $B_j\in\BB_j$.

\noindent
(iii) One has a semiorthogonal decomposition
$$D^b_G(X)=\lan \pi^*D^b(Y), \BB_1,\ldots,\BB_r \ran.$$
\end{thm}

\Pf . (i) For $F_1,F_2\in D^b(X)$ we have
$$\Hom_G(\pi^*F_1,\pi^*F_2)\simeq\Hom(F_1,\pi_*\pi^*F_2)^G
\simeq\Hom(F_1,F_2\ot (\pi_*\OO_X)^G)\simeq\Hom(F_1,F_2),$$
since $(\pi_*\OO_X)^G\simeq\OO_Y$.

\noindent
(ii) Let us first prove that the collection of $G$-sheaves on $X$
\begin{equation}\label{exc-col2}
(\OO_{D_i}, \zeta(i)_{D_i}[1],\ldots,\zeta(i)^{m_i-2}_{D_i}[m_i-2])
\end{equation}
is exceptional. Indeed, it is clear that there are no $G$-morphisms between $\zeta(i)^a_{D_i}$
and $\zeta(i)^b_{D_i}$ for $a\not\equiv b\mod(m_i)$.
Also, by Serre duality, for $a,b\in\Z/m_i\Z$ we have
$$\Ext^1_G(\zeta(i)^a_{D_i},\zeta(i)^b_{D_i})^*\simeq
\Hom_G(\zeta(i)^b_{D_i},\zeta(i)^{a+1}_{D_i})\simeq\Hom_{G_i}(\zeta(i)^b,\zeta(i)^{a+1}).$$
The latter space is nonzero only for $b=a+1$. This proves that \eqref{exc-col2} is
exceptional.
Now using the exact sequences 
$$0\to\zeta(i)^a_{D_i}\to\OO_{(a+1)D_i}\to\OO_{aD_i}\to 0$$
for $a=1,\ldots,m_i-2$, one can easily show that making a sequence of mutations
in \eqref{exc-col2} one gets the sequence \eqref{exc-col1}. 
%More precisely, \eqref{exc-col1} is the left (right?) dual exceptional collection to \eqref{exc-col2}. 
The fact that $\BB_i$ and $\BB_j$ are mutually orthogonal follows
from disjointness of $D_i$ and $D_j$.

\noindent
(iii) Since the subcategory $\lan\BB_1,\ldots,\BB_r\ran$ is admissible
it is enough to prove that $\pi^*D^b(Y)$ coincides with its right orthogonal.
Since $\BB_i$ is also generated by the exceptional collection \eqref{exc-col2},
the condition $\Hom_G(\BB_i,F)=0$ for $F\in D^b_G(X)$ is equivalent
$$\Hom_G^*(\zeta(i)^a_{D_i},F)=0\text{ for }a=0,\ldots,m_i-2.$$
Using Serre duality we can rewrite this as
$$\Hom_G(F,\zeta(i)^a_{D_i})=0\text{ for }a=1,\ldots,m_i-1.$$
Equivalently, $F|_{p_i}$ should have a trivial $G_i$-action.
By the main theorem of \cite{T} this implies that $F\in\pi^*D^b(Y)$.
\ed

In the case when  $Y=\P^1$ the semiorthogonal decomposition of Theorem \ref{semiort-thm} gives rise to a full exceptional collection in $D^b_G(X)$.

\begin{cor} \label{exceptional-cor} Assume $X/G\simeq\P^1$. Then for every $n\in\Z$ we have 
the following full exceptional collection in $D^b_G(X)$:
$$(\pi^*\OO_{\P^1}(n),\pi^*\OO_{\P^1}(n+1),\OO_{(m_1-1)D_1},\ldots,\OO_{D_1},\ldots,
\OO_{(m_r-1)D_r},\ldots,\OO_{D_r}).$$
In particular, $K^0(D^b_G(X))\simeq\Z^{2+\sum_{i=1}^r(m_i-1)}$.
\end{cor}

\begin{defi}
For a collection of $r$ distinct points $\ov{\la}=(\la_1,\ldots,\la_r)$ on $\P^1(k)$ and a sequence of weights
$\ov{m}=(m_1,\ldots,m_r)$ let us define the algebra
$\La(\ov{\la},\ov{m})$ as the path algebra of a quiver $Q_{\ov{m}}$ modulo 
relations $I(\ov{\la})$, where

\noindent
(i) $Q_{\ov{m}}$ has $2+\sum_{i=1}^r(m_i-1)$ vertices named $u$, $v$, and $w_1^1,\ldots,w_1^{m_1-1},\ldots,w_r^1,\ldots,w_r^{m_r-1}$;

\noindent
(ii) $Q_{\ov{m}}$ has $2+\sum_{i=1}^r(m_i-1)$ arrows: $2$ arrows $u\stackrel{x_0,x_1}{\to}v$,
and chains of arrows
$$v\stackrel{e_i}{\to}w_i^{m_i-1}\to w_i^{m_i-2}\to\ldots\to w_i^1$$
for every $i=1,\ldots,r$; 

\noindent
(iii) $I(\ov{\la})$ is generated by $r$ quadratic relations: $L_i\cdot e_i=0$, $i=1,\ldots,m$,
where $L_i\sub kx_0\oplus kx_1$ is the line corresponding to $\la_i\in\P^1=\P(kx_0\oplus kx_1).$
\end{defi}

It is easy to see that the endomorphism algebra of the exceptional collection constructed
in Corollary \ref{exceptional-cor} is isomorphic to $\La(\ov{\la},\ov{m})$,
where $D_i=\pi^{-1}(\la_i)$.
Hence, we obtain the following description of
the derived category of $G$-sheaves on a $G$-covering of $\P^1$.
(where for a finite-dimensional
algebra $A$ we denote by $\mod-A$ the category of finite-dimensional right $A$-modules).

\begin{cor}\label{quiver-cor} Let $\pi:X\to\P^1$ be a ramified Galois covering, 
where $X$ is a smooth curve, with Galois group $G$.
Let $\ov{\la}=(\la_1,\ldots,\la_r)\sub\P^1$ be the set of ramification points of $\pi$ and let
$\ov{m}=(m_1,\ldots,m_r)$ be multiplicities of the corresponding fibers. 
Then for every $n\in\Z$ one has an exact equivalence of triangulated categories
$$\Phi_n:D^b_G(X)\wt{\to} D^b(\mod-\La(\ov{\la},\ov{m})):
F\to R\Hom_G(V_n,F),$$
where 
$$V_n=\pi^*\OO_{\P^1}(n)\oplus\pi^*\OO_{\P^1}(n+1)\oplus\bigoplus_{1\le i\le r, 1\le j<m_i}
\OO_{jD_i}$$
\end{cor}

\begin{rem} Another natural full exceptional collection in $D^b_G(X)$ is
\begin{equation}\label{exc-line}
(\OO_X, \OO_X(D_1),\ldots,\OO_X((m_1-1)D_1),\ldots,\OO_X(D_r),\ldots,\OO_X((m_r-1)D_r),
\pi^*\OO_{\P^1}(1)).
\end{equation}
It is obtained from the collection of Corollary \ref{exceptional-cor} for $n=1$ by making the left mutation through $\pi^*\OO_{\P^1}(1)$ of the part of the collection following this object.
 \end{rem} 

A right module $M$ over
$\La(\ov{\la},\ov{m})$ 
can be viewed as a representation of the quiver
$Q_{\ov{m}}^{op}$ in which the relations $I(\ov{\la})^{op}$ are satisfied. Thus, $M$
is given by a collection of vector spaces $(U,V,W_i^j)$,
$i=1,\ldots,r$, $j=1,\ldots,m_i-1$,
equipped with linear maps
$$W_i^1\to\ldots\to W_i^{m_i-1}\to V, \ i=1,\ldots,r,$$ 
and $x_0,x_1:V\to U$
satisfying the relations $I(\ov{\la})^{op}$. Let us define the following additive functions of $M$:
$$\deg_n(M)=|G|\cdot\left(n\dim U-(n-1)\dim V-\sum_{i=1}^r\sum_{j=1}^{m_i-1}\frac{\dim W_i^j}{m_i}
\right),$$
$$\rk(M)=\dim U-\dim V.$$
We extend these functions to additive functions on $D^b(\La(\ov{\la},\ov{m}))$.

\begin{lem}\label{deg-rk-lem} 
In the situation of Corollary \ref{quiver-cor} one has 
$\deg_n(\Phi_n(F))=\deg(F)$ and $\rk(\Phi_n(F))=\rk(F)$
for every $F\in D^b_G(X)$. 
\end{lem}

\Pf . Let $M=\Phi_n(F)=R\Hom_G(V_n,F)$ for $F\in D^b_G(X)$.
Then we have
$$\dim U=\chi_G(\pi^*\OO_{\P^1}(n),F),\ \dim V=\chi_G(\pi^*\OO_{\P^1}(n+1),F),\
\dim W_i^j=\chi_G(\OO_{jD_i},F),$$
and our task is to express $\rk(F)$ and $\deg(F)$ in terms of these numbers. 
To compute the rank we can use the equality
$$\rk(F)=-\chi_G(\pi^*\OO_q,F),$$
where $q$ is a generic point of $\P^1$. Since
$[\OO_q]=[\OO_{\P^1}(n+1)]-[\OO_{\P^1}(n)]$, this immediately implies
the required formula
$$\rk(F)=\chi_G(\pi^*\OO_{\P^1}(n),F)-\chi_G(\pi^*\OO_{\P^1}(n+1),F).$$
The formula for $\deg(F)$ should have form
$$\deg(F)=a\chi_G(\pi^*\OO_{\P^1}(n),F)+b\chi_G(\pi^*\OO_{\P^1}(n+1),F)+
\sum_{i=1}^r\sum_{j=1}^{m_i-1}c_i^j\chi_G(\OO_{jD_i},F)$$
for some constants $a$, $b$ and $c_i^j$.
The constants are determined by substituting in this formula the elements of the dual basis of $K_0(D^b_G(X))$:
$$([\pi^*\OO_{\P^1}(n)],-[\pi^*\OO_{\P^1}(n-1)], -[\zeta(1)_{D_1}],\ldots,-[\zeta(1)^{m_1-1}_{D_1}],
\ldots,-[\zeta(r)_{D_r}],\ldots,-[\zeta(r)^{m_r-1}_{D_r}]).$$
\ed

\subsection{Elliptic Galois coverings of $\P^1$}\label{ell-sec}

Now let us specialize to the case of a Galois covering $\pi:E\to\P^1$, where $E=\C/(\Z\tau+\Z)$ is an elliptic curve (so $k=\C$). More precisely, we are interested in the following four cases in which
$G$ is a cyclic subgroup in $\C^*$ acting on $E$ in the natural way.

\noindent
(i) $E$ is arbitrary and $G=\Z/2\Z$. 
The corresponding double covering
$\pi:E\to\P^1$ is given by the Weierstrass $\wp$-function and is ramified exactly over $4$-points
$$\{\la_1,\la_2,\la_3,\la_4\}=\{\infty, \wp(\frac{1}{2}),\wp(\frac{\tau}{2}),\wp(\frac{1+\tau}{2})\}.$$

\noindent
(ii) $E=\C/L_{tr}$, where $L_{tr}=\Z\frac{1+i\sqrt{3}}{2}+\Z$, and $G=\Z/3\Z$.
In this case $\pi:E\to\P^1$ is given by $\wp'(z)$. Note that $E^G$ consists of $3$ points:
$0\mod L_{tr}$ and $\pm\frac{3+i\sqrt{3}}{6}\mod L_{tr}$. Hence, $\pi$ is totally ramified over $3$ points.

\noindent
(iii) $E=\C/L_{sq}$, where $L_{sq}=\Z i+\Z$, and $G=\Z/4\Z$.
 In this case $\pi:E\to\P^1$ is given by $\wp(z)^2$. 
We have two points whose stabilizer subgroup is $\Z/2\Z$, namely, $\frac{1}{2}\mod L_{sq}$ and 
$\frac{i}{2}\mod L_{sq}$ (they get exchanged by the generator of $\Z/4\Z$). The two points in $E^G$
are $0\mod L_{sq}$ and $\frac{i+1}{2}\mod L_{sq}$. Hence, $\pi$ is ramified over $3$ points and
the corresponding multiplicities are $(2,4,4)$.

\noindent
(iv) $E=\C/L_{tr}$ (same curve as in (ii)) and $G=\Z/6\Z$. In this case
$\pi:E\to\P^1$ is given by $\wp'(z)^2$. There are $3$ points whose stabilizer subgroup is
$\Z/2\Z$, namely, all nontrivial points of order $2$ on $E$ (they form one $G$-orbit).
There is also a $G$-orbit consisting of two points $\pm\frac{3+i\sqrt{3}}{6}\mod L_{tr}$
with stabilizer subgroup $\Z/3\Z$. Finally, $0\mod L_{tr}$ is the only point in $E^G$.
Therefore, $\pi$ is ramified over $3$ points with multiplicities $(2,3,6)$.

From the above description of the ramification data
and from Corollary \ref{exceptional-cor} we get 

\begin{cor}\label{ell-cor}
One has $K_0(D^b_G(E))\simeq\Z^r$, where $r=6,8,9$ or $10$ 
in the cases (i)-(iv), respectively.
\end{cor}

\subsection{Galois coverings of $\P^1$ and weighted projective curves}

The results of this section are not used in the rest of the paper. Its purpose is to
explain the relation between $G$-sheaves on ramified Galois coverings
of $\P^1$ and coherent sheaves on weighted projective curves introduced in \cite{GL}. 
This relation is known to experts, however, our proof seems to be new.

Let us recall the definition of weighted curves 
\footnote{These curves are also called
weighted projective lines, however, they should not be confused with one-dimensional weighted projective spaces.}
$C(\ov{m},\ov{\la})$ of \cite{GL} 
associated with a sequence of positive integers $\ov{m}=(m_1,\ldots,m_r)$ and a
sequence of points $\ov{\la}=(\la_1,\ldots,\la_r)$ in $\P^1(k)$. 
Let $Z(\ov{m})$ be the rank one abelian group with generators $e_1,\ldots,e_r$ and relations
$m_1e_1=\ldots=m_re_r$. Let us also choose for every $i=1,\ldots,r$ a nonzero section
$s_i\in H^0(\P^1,\OO_{\P^1}(1))$ such that $s_i(\la_i)=0$.
Consider the algebra 
$$S(\ov{m},\ov{\la})=k[x_1,\ldots,x_r]/I(\ov{m},\ov{\la}),$$
where the ideal $I(\ov{m},\ov{\la})$ is generated by all polynomials of the form
$a_1x_1^{m_1}+\ldots a_r x_r^{m_r}$ such that
$\sum_{i=1}^r a_i s_i=0$. 
Let $Z(\ov{m})_+\sub Z(\ov{m})$ be the positive submonoid generated by $e_1,\ldots,e_r$.
Note that the algebra $S(\ov{m},\ov{\la})$ has a natural $Z(\ov{m})_+$-grading, where
$\deg(x_i)=e_i$.
The category $\coh(C(\ov{m},\ov{\la}))$ of coherent sheaves on $C(\ov{m},\ov{\la})$ can be
defined as the quotient-category of the category of finitely generated
$Z(\ov{m})_+$-graded $S(\ov{m},\ov{\la})$-modules by the subcategory of finite length modules.

Now assume we are given a ramified Galois covering $\pi:X\to\P^1$ with Galois group
$G$, where $X$ is a smooth connected curve. Define the associated data $(D_i)$, $\ov{m}$ and $\ov{\la}$ as in  the previous section. Let $\Pic_G(X)$ be the group of $G$-equivariant line bundles 
up to $G$-isomorphism. Let us consider the algebra
$$S(X,G):=\oplus_{[L]\in\Pic_G(X)}H^0(X,L)^G.$$

\begin{thm}\label{isom-alg-thm} 
One has an isomorphism of algebras
$$S(\ov{m},\ov{\la})\simeq S(X,G),$$
compatible with gradings via an isomorphism $Z(\ov{m})\simeq \Pic_G(X)$.
\end{thm}

\Pf . We claim that there is a natural homomorphism $S(\ov{m},\ov{\la})\to S(X,G)$ that sends
$x_i$ to a nonzero section $f_i$ of $H^0(X,\OO_X(D_i))^G$. Indeed,
note that we have a natural isomorphism $\OO_X(m_iD_i)\simeq\pi^*\OO_{\P^1}(1)$
compatible with the action of $G$, and hence the induced isomorphism
$$H^0(\OO_X(m_iD_i))^G\simeq H^0(\P^1,\OO_{\P^1}(1)).$$
Let us rescale $f_i$ in such a way that $f_i^{m_i}$ corresponds to $s_i\in H^0(\P^1,\OO_{\P^1})$
under this isomorphism. Then $f_i^{m_i}$ will satisfy the same linear relations as $s_i$, hence
we get a homomorphism $\a:S(\ov{m},\ov{\la})\to S(X,G)$. Note that $\a$ is compatible
with gradings via the homomorphism $Z(\ov{m})\to\Pic_G(X)$ sending $x_i$ to the class of $D_i$. Let us check that $\a$ is surjective. Assume we are given $L\in\Pic_G(X)$ and a nonzero $G$-invariant
section $f$ of $L$. If the divisor of zeroes of $f$ contains $D_i$ for some $i$ then
$f$ is divisible by $f_i$ in the algebra $S(X,G)$, so we can assume that the divisor of
$f$ is disjoint from all special fibers. Therefore, $L\simeq\pi^*\OO_{\P^1}(n)$ and
$f$ corresponds to a section of $\OO_{\P^1}(n)$ on $\P^1$. Note that $r\ge 2$
since $X$ is connected. Therefore, every section of $\OO_{\P^1}(n)$ can be expressed as
a polynomial of $s_1,\ldots,s_r$. Hence, $f$ belongs to the image of $\a$.
Injectivity of $\a$ follows easily from Proposition 1.3 of \cite{GL}.
\ed

Using Theorem \ref{isom-alg-thm} we can derive the following
equivalence between the categories of sheaves. 

\begin{thm}\label{Serre-thm} 
In the above situation one has an equivalence of categories
$$\Coh_G(X)\simeq\Coh(C(\ov{m},\ov{\la})).$$
\end{thm}

\Pf . This follows from Theorem \ref{isom-alg-thm} by a version of Serre's theorem. 
The only nontrivial
fact one has to use is that for every $G$-sheaf $F$ on $X$ there exists a surjection
of $G$-sheaves $\oplus_{i=1}^n L_i\to F$, where $L_i$ are equivariant $G$-bundles.
Since every $G$-sheaf $F$ is covered
by a $G$-bundle of the form $H^0(X,F\otimes L)\otimes L^{-1}$ for sufficiently ample
$G$-equivariant line bundle $L$, it suffices to consider the case when $F$ is locally free.
Assume first that the action of $G_i$ on the fiber $F|_{p_i}$ is trivial for all $i=1,\ldots,r$.
Then $F$ is $G$-isomorphic to the pull-back of a vector bundle on $\P^1$.
In this case the assertion is
clear since all vector bundles on $\P^1$ are direct sums of line bundles.
We are going to reduce to this case using elementary transformations along $D_i$'s.
Namely, let us decompose a representation of $G_i$ on the fiber $F|_{p_i}$ into
the direct sum of characters of $G_i$:
$$F|_{p_i}\simeq \oplus_{j=0}^{m_i-1}V_j\otimes\zeta(i)^j$$
with some multiplicity spaces $V_j$. Note that if we define the $G$-bundle $F'$ by the short exact triple
$$0\to F'\to F\to V_j\otimes\zeta(i)^j_{D_i}\to 0$$
then we have an exact sequence of $G_i$-modules
$$0\to V_j\otimes \zeta(i)^j\otimes\Tor_1(\OO_{D_i},\OO_{p_i})\to F'|_{p_i}\to F_{p_i}\to V_j\otimes\zeta(i)^j\to 0.$$
But $\Tor_1(\OO_{D_i},\OO_{p_i})\simeq\om_X|_{p_i}\simeq\zeta(i)$, hence,
$$F'|_{p_i}\simeq \left(\oplus_{j'\neq j} V_{j'}\otimes\zeta(i)^{j'}\right)\oplus V_j\otimes\zeta(i)^{j+1}.$$
It is clear that using a sequence of transformations of this form we can pass from $F$ to a
vector bundle for which all fibers $F|_{p_i}$ have trivial $G_i$-action. It remains to check
that if our claim holds for $F'$ (i.e., there exists a $G$-surjection from a direct sum of $G$-equivariant
line bundles to $F'$) then the same is true for $F$. To this end we observe that for every $n\in\Z$
there exists a surjection of $G$-sheaves 
$$\om_X^j\otimes\pi^*\OO_{\P^1}(n)\to\zeta(i)^j_{D_i}.$$
If $n$ is sufficiently negative then this map lifts to a morphism
$\om_X^j\otimes\pi^*\OO_{\P^1}(n)\to F.$
Thus, from a surjection $\oplus_{i=1}^n L_i\to F'$ we obtain a surjection
of the form
$$\oplus_{i=1}^n L_i\oplus \om_X^j\otimes\pi^*\OO_{\P^1}(n)\to F.$$
\ed

In the case when $G=\Z/2\Z$ and
$X$ is an elliptic curve the above equivalence is considered in Example 5.8 of \cite{GL}.

Note that the tilting bundle on $C(\ov{m},\ov{\la})$ constructed
in \cite{GL} corresponds to the exceptional collection \eqref{exc-line}.

\section{Holomorphic bundles on toric orbifolds and derived categories of $G$-sheaves}

\subsection{Remarks on $B_{\th}$-modules and holomorphic bundles}\label{remarks-sec}

It is clear that a $B_{\th}$-module is finitely generated iff it is finitely generated as an $A_{\th}$-module.
We claim that projectivity also can be checked over $A_{\th}$.

\begin{lem}\label{proj-lem} 
Let $M$ be a right $B_{\th}$-module. Then $M$ is projective as a $B_{\th}$-module
iff it is projective as an $A_{\th}$-module.
\end{lem}

\Pf . The ``only if" part is clear. Let $M$ be a right $B_{\th}$-module, projective over $A_{\th}$.
Then we have a natural surjection of $B_{\th}$-modules
$p:M\otimes_{A_{\th}}B_{\th}\to M$ given by the action of
$B_{\th}$. On the other hand, it is easy to check that the map
$$s:M\to M\otimes_{A_{\th}}B_{\th}:m\mapsto \frac{1}{|G|}\sum_{g\in G} mg\otimes g^{-1}$$
commutes with the right action of $B_{\th}$. Since $p\circ s=\id_M$, we derive that
$M$ is a direct summand in the projective $B_{\th}$-module $M\otimes_{A_{\th}}B_{\th}$.
Hence, $M$ itself is a projective $B_{\th}$-module.
\ed

Thus, we can identify holomorphic bundles on $T_{\th,\tau}/G$ with $G$-equivariant
holomorphic bundles on $T_{\th,\tau}$. Here is a more precise statement.
Let us define an automorphism $g^*$ of the category $\Hol(T_{\th,\tau})$ by
setting $g^*(P,\ov{\nabla})=(P^{g},\vareps(g)\ov{\nabla})$, where $P^{g}=P$ as a vector space
but the $A_{\th}$-module structure is changed by the automorphism $g$ of $A_{\th}$. The
fact that we again obtain a holomorphic bundle on $T_{\th,\tau}$ follows from \eqref{gdelta}.

\begin{lem}\label{eq-hol-bun-lem} 
The category $\Hol(T_{\th,\tau}/G)$ is equivalent to
the category of $G$-equivariant objects of $\Hol(T_{\th,\tau})$.
\end{lem}

\Pf . By Lemma \ref{proj-lem} a holomorphic bundle on $T_{\th,\tau}/G$ is given by a
finitely generated projective right $A_{\th}$-modules $P$ equipped with a holomorphic
structure $\ov{\nabla}$ and an action of $G$ such that 
$$g(f\cdot a)=g(f)\cdot g(a),\ g\circ\ov{\nabla}=\vareps(g)\ov{\nabla}\circ g,$$
where $g\in G$, $f\in P$, $a\in A_{\th}$. This immediately implies the assertion.
\ed

\begin{prop}\label{hol-str-prop}
Every finitely generated projective $B_{\th}$-module admits a holomorphic
structure.
\end{prop}

\Pf . Let $P$ be such a module. Considering $P$ as an $A_{\th}$-module we can equip
it with a holomorphic structure $\ov{\nabla}$ making it into a holomorphic bundle on $T_{\th}$
(because $P$ is a direct sum of basic modules and every basic module admits a 
standard holomorphic structure, see \cite{PS}). Now replace $\ov{\nabla}$ with
$$\frac{1}{|G|}\sum_{g\in G}\eps(g)^{-1} g\ov{\nabla}g^{-1}.$$ 
This new structure is compatible with the
action of $G$,
so that $P$ becomes a holomorphic bundle on $T_{\th,\tau}/G$.
\ed

\subsection{Generalities on torsion theory}

Recall (see \cite{HRS}) that a {\it torsion pair} in an exact category $\CC$ is a pair
of full subcategories $(\TT,\FF)$ in $\CC$ such that $\Hom(T,F)=0$ for every 
$T\in\TT$, $F\in\FF$, and every object $C\in\CC$ fits into a short exact triple
$$0\to T\to C\to F\to 0$$
with $T\in\TT$ and $F\in FF$. 
Note that if $(\TT,\FF)$ is a torsion pair then $\FF$ (resp., $\TT$) coincides with the right (resp., left)
orthogonal of $\TT$, i.e. with the full subcategory of objects $X$ such that
$\Hom(T,X)=0$ for all $T\in\TT$ (resp., $\Hom(X,F)=0$ for all $F\in\FF$).
In particular $\TT$ and $\FF$ are stable under extensions and passing to a direct summand.

It will be convenient for us to introduce a slight generalization of the notion of a torsion pair.
Given a collection of full subcategories $(\CC_1,\ldots,\CC_n)$ in an exact category $\CC$,
let us denote by $[\CC_1,\ldots,\CC_n]$ the full subcategory in $\CC$ consisting of objects
$C$ admitting an admissible filtration
$0=F_0C\sub F_1C\sub\ldots\sub F_nC=C$ such that $F_iC/F_{i-1}C\in\CC_i$ for
$i=1,\ldots,n$.

\begin{defi} A {\it torsion $n$-tuple} in an exact category $\CC$ is a collection of
full subcategories $(\CC_1,\ldots,\CC_n)$ such that $\Hom(C_i,C_j)=0$ whenever
$C_i\in\CC_i$, $C_j\in\CC_j$, $i<j$, and $[\CC_1,\ldots,\CC_n]=\CC$. 
\end{defi}

Sometimes we will write the condition of absence of nontrivial morphisms in the above
definition as $\Hom(\CC_i,\CC_j)=0$ for $i<j$.
For $n=2$ we recover the notion of a torsion pair. Moreover, it is clear that if
$(\CC_1,\ldots,\CC_n)$ is a torsion $n$-tuple then for every $i$ the pair
$$([\CC_1,\ldots,\CC_i],[\CC_{i+1},\ldots,\CC_n])$$
is a torsion pair.
Note that the subcategories $\CC_i$ in this definition are automatically stable under extensions.
The main reason for introducing torsion $n$-tuples is because it is possible to substitute one
such torsion tuple into another. Namely, if $(\CC_1,\ldots,\CC_n)$ is a torsion $n$-tuple in $\CC$,
and $(\CC_{i,1},\ldots,\CC_{i,m})$ is a torsion $m$-tuple in $\CC_i$ then
$$(\CC_1,\ldots,\CC_{i-1},\CC_{i,1},\ldots,\CC_{i,m},\CC_{i+1},\ldots,\CC_n)$$
is a torsion $(n+m-1)$-tuple in $\CC$.

If $\CC$ is abelian then a torsion pair $(\TT,\FF)$ defines a nondegenerate $t$-structure on the derived category $D^b(\CC)$ with the heart
$$\CC^p:=\{K\in D^b(\CC):\  H^i(K)=0\text{ for }i\neq 0,-1, H^0(K)\in\TT, H^{-1}(K)\in\FF\}$$
(see \cite{HRS}). In other words,
$$\CC^p=[\FF[1],\TT],$$
where for a pair of full subcategories $\CC_1,\CC_2$ in a triangulated category $\DD$
we denote by $[\CC_1,\CC_2]$ the full subcategory in $\DD$ consisting of objects $K$
that fit into an exact triangle
$$C_1\to K\to C_2\to C_1[1]$$
with $C_1\in\CC_1$, $C_2\in\CC_2$.
The process of passing from $\CC$ to $\CC^p$ is called {\it tilting} (also, we
will call $\CC^p$ a {\it tilt} of $\CC$). Note that $(\FF[1],\TT)$ is a torsion pair in $\CC^p$
and applying tilting to this pair we pass back to $\CC$.
If $(\CC_1,\ldots,\CC_n)$ is a torsion $n$-tuple in $\CC$ then we set
$$[\CC_{i+1}[1],\ldots,\CC_n[1],\CC_1,\ldots,\CC_i]:=
[[\CC_{i+1},\ldots,\CC_n][1],[\CC_1,\ldots,\CC_i]]\sub\D^b(\CC),$$
where $([\CC_1,\ldots,\CC_i],[\CC_{i+1},\ldots,\CC_n])$ is the corresponding
torsion pair in $\CC$.

The main example relevant for us is the torsion pair 
$(\Coh_{>\th}(X),\Coh_{<\th}(X))$ in the category $\Coh(X)$ of coherent sheaves
on a smooth projective curve $X$, associated with an irrational number $\th$. 
Namely,
$\Coh^{<\th}(X)\sub\Coh(X)$ (resp., $\Coh^{>\th}(X)\sub\Coh(X)$) consists
of all coherent  sheaves $F$ on $X$ 
such that all subsequent quotients in the 
Harder-Narasimhan filitration of $F$ have slope $<\th$ (resp., $>\th$), where we consider
torsion sheaves as having slope $+\infty$.
Note that these torsion pairs arise in connection with stability structures on $D^b(X)$
(see \cite{Bridge}). 
%In the case of rational $\th$ one can consider two torsion pairs
%$(\Coh_{\ge \th}(X),\Coh_{<\th})(X)$ and $(\Coh_{>\th}(X),\Coh_{\le \th}(X))$
%defined in a similar fashion.

\subsection{Fourier-Mukai transform for noncommutative two-tori}

Let 
$$\CC^{\th}(E)=[\Coh^{<\th}(E)[1],\Coh^{>\th}(E)]\sub D^b(E)$$
be the tilt of the category of coherent sheaves on the elliptic curve $E=\C/(\Z\tau+\Z)$
associated with $\th$. 
We know from \cite{PS},\cite{P} that $\Hol(T_{\th,\tau})$ is equivalent to $\CC^{\th}(E)$.
In this section we will show that the construction of this equivalence can be adjusted to
be compatible with the action of a finite group $G$.

Recall that the equivalence is given by a version of the Fourier-Mukai transform (see 
\cite{PS}, Section 3.3) . With a holomorphic vector bundle $(P,\ov{\nabla})$ on $T_{\th,\tau}$
this transform associates the complex $\SS(P,\ov{\nabla})$ of $\OO$-modules on $E$
of the form $d:P_E\to P_E$ concentrated in degrees $[-1,0]$, where $P_E$ is obtained by descending the sheaf
of holomorphic $E$-valued functions over $\C$ using an action of $\Z^2$ of the form
$$\rho_v(f)(z)=\exp(\pi i\th c_v(z)) f(z+v)U_v, \ v\in\Z^2,$$
and the differential $d$ is induced by the operator 
$$f(z)\mapsto\ov{\nabla}(f(z))+2\pi iz f(z)$$.
Here $(c_v(z))$ is a collection of holomorphic functions on $\C$ numbered by $\Z^2$ satisfying
the condition
\begin{equation}\label{Fourier-cocycle}
c_{v_1}(z)+c_{v_2}(z+v_1)-c_{v_1+v_2}(z)=\det(v_1,v_2).
\end{equation}
Note that in \cite{PS} we made one possible choice of $(c_v(z))$, however, it is not the only choice.
In fact, one can easily see that the equivalences corresponding to different choices of $(c_v(z))$ differ
by tensoring with a holomorphic line bundle on $E$.
One of possible solutions of \eqref{Fourier-cocycle} is
$$c^0_{m\tau+n}(z)=-2mz-m(m\tau+n).$$
It follows from Proposition 3.7 of \cite{PS} that $\SS$ is an equivalence of $\Hol(T_{\th,\tau})$ with $\CC^{\th}(E)$ (note that the definition above differs from that of \cite{PS} by the shift of degree).

Now let us assume that a finite group $G$ acts on the elliptic
curve $E$ be automorphisms (preserving $0$). This means that $G$ is a subgroup in $\C^*$
and multiplication by elements of $G$ preserves the lattice $\Z\tau+\Z\sub\C$. Identifying
$\Z^2$ with $\Z\tau+\Z$ by $(m,n)\mapsto m\tau+n$ we can view $G$ also as a subgroup in
$\SL_2(\Z)$. One can immediately check that the corresponding action of $G$ on $A_{\th}$ satisfies
\eqref{gdelta} with $\vareps(g)=g^{-1}\in\C^*$.
Hence, for every $g\in G$ we have the corresponding automorphism $g^*$
of the category $\Hol(T_{\th,\tau})$ (see section \ref{remarks-sec}). 
Let us make a $G$-invariant choice of $(c_v(z))$ by setting
$$c_v(z)=\frac{1}{|G|}\sum_{g\in G} c^0_{gv}{gz},$$
so that $c_{gv}(gz)=c_v(z)$ for all $g\in G$. Then the resulting Fourier-Mukai transform $\SS$ is compatible with the action of $G$ in the standard way.

\begin{prop}\label{G-Fourier-prop} 
With the above choice of $(c_v(z))$ one has natural isomorphisms of functors
$$\SS\circ g^*\simeq (g^{-1})^*\circ\SS$$
from $\Hol(T_{\th,\tau})$ to $\CC^{\th}(E)$, where $g\in G$.
\end{prop}

\Pf . By definition $g^*(P,\ov{\nabla})=(P^g,\vareps(g)\ov{\nabla})$.
Hence, $\SS g^*(P,\ov{\nabla})=[d_1:P^g_E\to P^g_E)]$ where $P^g_E$ is obtained 
from the action of $\Z^2$ on $P_{\C}$ given by
$$f(z)\mapsto \exp(\pi i \th c_v(z))f(z+v)U_{gv},$$
and the differential $d_1$ is induced by the operator 
$$f(z)\mapsto \vareps(g)\ov{\nabla}(f(z))+2\pi i zf(z).$$
On the other hand, $(g^{-1})^*\SS(P,\ov{\nabla})$ is given by the complex
$[d_2:(g^{-1})^*P_E\to (g^{-1})^*P_E]$, where $(g^{-1})^*P_E$ is obtained from the
action of $\Z^2$ on $P_{\C}$ given by
$$f(z)\mapsto \exp(\pi i \th c_v(gz))f(z+g^{-1}v)U_v,$$
and $d_2$ is induced by the operator
$$f(z)\mapsto \ov{\nabla}(f(z))+2\pi i gz f(z),$$
where we view $g$ as an element of $\C^*$. Making a change of variables $v\mapsto gv$
we can identify two $\Z^2$-actions above, and hence we can identify with $P^g_E$
with $(g^{-1})^*P_E$. Since $\vareps(g)=g^{-1}$, under this identification
$d_1=\vareps(g)d_2$, so we get the required isomorphism. 
\ed

\subsection{Proof of Theorem \ref{mainthm}}

From Lemma \ref{eq-hol-bun-lem} and Proposition \ref{G-Fourier-prop} we obtain that
the category $\Hol(T_{\th,\tau}/G)$ is equivalent to
the category of $G$-equivariant objects of $\CC^{\th}(E)$.

Note that the Harder-Narasimhan filtration of a $G$-sheaf
is stable under the action of $G$ and hence, can be considered as a filtration in
$\Coh_G(E)$. Therefore, we can define a torsion theory
$((\Coh^{>\th}_G(E),\Coh^{<\th}_G(E))$ in $\Coh_G(E)$, where
$\Coh^*_G(E)$ consists of $G$-sheaves $F$ such that
after forgetting the $G$-structure we have $F\in\Coh^*(E)$. 
Let 
$$\CC^{\th}_G(E)=[\Coh^{<\th}_G(E)[1],\Coh^{>\th}_G(E)]\sub D_G^b(E)$$ 
be the corresponding tilted abelian subcategory.
By Lemma \ref{eq-sh-lem} the category $\CC^{\th}_G(E)$ is equivalent
to the category of $G$-equivariant objects on $\CC^{\th}(E)$, and hence, to
$\Hol(T_{\th,\tau}/G)$. 

Let us show that
$D^b_G(E)$ is equivalent to $D^b(\CC^{\th}_G(E))$.
By Proposition 5.4.3 of \cite{BVdB} 
it is enough to check that that our torsion pair in $\Coh_G(E)$ is {\it cotilting}, i.e.,
for every $G$-sheaf $F$ on $E$ there exists
a $G$-equivariant vector bundle $V\in\Coh^{<\th}(E)$ and a $G$-equivariant
surjection $V\to F$. However, this is clear since for every $G$-sheaf $F$ there is
a surjection
$$H^0(X,F\otimes\pi^*\OO_{\P^1}(n))\otimes \pi^*\OO_{\P^1}(-n)\to F,$$
where $n$ is large enough (and the space of global sections is equipped
with the natural $G$-action).

Thus, we showed that $\Hol(T_{\th,\tau})/G$ is abelian and its derived category is equivalent
to $D^b_G(E)$. It remains to apply Corollary \ref{quiver-cor} to the covering
$\pi: E\to E/G\simeq\P^1$.
The statement about the number of vertices follows from the explicit description of these
coverings in section \ref{ell-sec}.
\ed

\subsection{Proof of Theorem \ref{K-thm}}

Using Theorem \ref{mainthm} and Corollary \ref{ell-cor}
 we obtain an isomorphism 
 $$K_0(\Hol(T_{\th,\tau}/G))\simeq K_0(D^b_G(E))\simeq\Z^r,$$ where
 $r=6,8,9,10$ for $G=\Z/m\Z$ with $m=2,3,4,6$, respectively.
 Now we observe that by Proposition \ref{hol-str-prop} the natural homomorphism
 $$K_0(\Hol(T_{\th,\tau}/G)\to K_0(B_{\th})$$
is surjective. To prove that this map is an isomorphism, it is enough to check
that the rank of $K_0(B_{\th})$ is at least $r$.
This was done in \cite{BEEK1},
 \cite{W-Fourier} and \cite{BW1} for $m=2$, $m=4$, and $m=3, 6$, respectively 
(by explicitly constructing $r$ elements in $K_0(B_{\th})$ and using unbounded
traces to check their linear independence). 
\ed

This result was known for $G=\Z/2\Z$ (see \cite{K}), however, with a different proof. 
For $G=\Z/4\Z$ and $G=\Z/6\Z$ it was known for $\th$ in a dense $G_{\delta}$-set
(see \cite{W-Ktheory} and \cite{BW2}). The case of $G=\Z/4\Z$ and general $\th$
was done by Lueck, Walters and Phillips (unpublished).

Note that from the above proof we also get the following 

\begin{cor} The natural homomorphism $K_0(\Hol(T_{\th,\tau}/G))\to K_0(B_{\th})$
is an isomorphism. Moreover, the positive cones are the same.
\end{cor}

\begin{rem}
In \cite{W} it was shown that for $G=\Z/2\Z$ the positive cone in $K_0(B_{\th})$ coincides with the
preimage of the positive cone in $K_0(A_{\th})$ under the natural homomorphism
$K_0(B_{\th})\to K_0(A_{\th})$ (in other words, it consists of all elements $x\in K_0(B_{\th})$
such that $\tr_*(x)>0$, where $\tr_*:K_0(B_{\th})\to\R$ is the homomorphism induced
by the trace). As was pointed to us by Chris Phillips, similar statement is also known to hold for
other groups $G$. Namely, it follows from the fact that the corresponding crossed products
are simple AH algebras with slow dimension growth and real rank zero (see Theorems 8.11 and
9.10 of \cite{Phillips}).
\end{rem}

\subsection{Tiltings associated with $\th$}

Let $\pi:X\to \P^1$ be a ramified Galois covering with the Galois group $G$, and let
$D^b_G(X)$ be the derived category
of $G$-sheaves on $X$. As in the above proof of Theorem \ref{mainthm} we can
define the torsion pair $(\Coh^{>\th}_G(X),\Coh^{<\th}_G(X))$ in 
$\Coh_G(X)$ associated with an irrational number $\th$. Our goal in this section is to describe the image of the corresponding tilted abelian subcategory
$$\CC^{\th}_G(X):=[\Coh^{<\th}_G(X)[1],\Coh^{>\th}_G(X)]\sub D^b_G(X)$$ 
under the equivalence
$\Phi_n$ of Corollary \ref{quiver-cor} (for suitable $n$).

Let us start by describing the torsion pair in $\Coh_G(X)$ giving rise to the $t$-structure on
$D^b_G(X)$ associated with $\Phi_n$. By definition, the heart $\MM_n$ of this $t$-structure
consists of objects $F$ such that $\Hom^i_G(V_n,F)=0$ for $i\neq 0$.

\begin{prop}\label{mod-tilt-prop} 
Let $\TT_0\sub\Coh_G(X)$ denote the full subcategory consisting of
all $G$-sheaves isomorphic to a direct sum of $G$-sheaves from the collection
$$(\OO_{(m_1-1)D_1},\ldots,\OO_{D_1},\ldots,\OO_{(m_r-1)D_r},\ldots,\OO_{D_r}).$$ 
Let also $\TT_1\sub\Coh_G(X)$ be the full subcategory
of torsion $G$-sheaves obtained by successive extensions from simple $G$-sheaves
of the form $\zeta(i)^a_{D_i}$, where $i=1,\ldots,r$, $a=1,\ldots,m_i-1$ (so $\OO_{D_i}$ are
not included). 
Then for every $n\in\Z$ we have a torsion quadruple
$$(\TT_0,\pi^*\Coh^{\ge n}(\P^1),\pi^*\Coh^{\le n-1}(\P^1),\TT_1)$$
in $\Coh_G(X)$. Furthermore, we have the equality of abelian subcategories in $D^b_G(X)$
\begin{equation}\label{mod-tilt-eq}
\MM_n=[\pi^*\Coh^{\le n-1}(\P^1)[1],\TT_1[1],\TT_0,\pi^*\Coh^{\ge n}(\P^1)].
\end{equation}
\end{prop}

\Pf . First, let us check that 
$$(\TT_0,\pi^*\Coh(\P^1),\TT_1)$$ 
is a torsion triple in $\Coh_G(X)$. The conditions
$$\Hom_G(\TT_0,\TT_1)=0,\ 
\Hom_G(\TT_0,\pi^*\Coh(\P^1))=0,\text{ and }\Hom(\pi^*\Coh(\P^1),\TT_1)=0$$
easily follow from the vanishings
$$\Hom_G(\OO_{bD_i},\zeta(i)^a_{D_i})=0,\
\Hom_G(\OO_{bD_i},\pi^*\OO_{\pi(p_i)})=0,\text{ and }
\Hom_G(\pi^*\Coh(\P^1),\zeta(i)^a_{D_i})=0,$$
where $a=1,\ldots,m_i-1$, $b=1,\ldots,m_i-1$.
Using elementary transformations along $D_i$'s as in the proof of Theorem
\ref{Serre-thm} one can easily see that every $G$-bundle on $X$ belongs
to $[\pi^*\Coh(\P^1),\TT_1]$. Now let $F$ be an indecomposable torsion
$G$-sheaf on $X$ supported on $D_i$. Then there exists a filtration
$$0=F_0\sub F_1\sub\ldots\sub F_n=F$$ 
by $G$-subsheaves such that 
$F_a/F_{a-1}\simeq \zeta(i)^{c-a}_{D_i}$
for $a=1,\ldots,n$ (where $c\in\Z/m_i\Z$). 
If $c-a\not\equiv 0\mod(m_i)$ for all $a=1,\ldots,n$ then $F\in\TT_1$.
Otherwise, let $a_1$ (resp., $a_2$) be the minimal (resp., maximal) $a$ 
such that $c-a\equiv 0\mod(m_i)$.
Then it is easy to see that
$$F_{a_1}\simeq\OO_{a_1D_i}\in\TT_0,\
F_{a_2}/F_{a_1}\in\pi^*\Coh(\P^1),\ F/F_{a_2}\in\TT_1,$$
and hence, $F\in [\TT_0,\pi^*\Coh(\P^1),\TT_1]$.

Substituting the torsion pair 
$$(\pi^*\Coh^{\ge n}(\P^1),\pi^*\Coh^{\le n-1}(\P^1))$$ 
into $\pi^*\Coh(\P^1)$
we obtain the required torsion quadruple. One can immediately check that 
each of the subcategories $\pi^*\Coh^{\le n-1}(\P^1)[1]$, $\TT_1[1]$, $\TT_0$ and
$\pi^*\Coh^{\ge n}(\P^1)$ belongs to $\MM_n$. Hence, the RHS of \eqref{mod-tilt-eq} is 
contained in $\MM_n$. Since both subcategories are hearts of nondegenerate $t$-structures,
this implies the required equality.
\ed

Note that $\pi^*\Coh^{\ge n}(\P^1)\sub\Coh^{\ge 2n}_G(X)$ and
$\pi^*\Coh_{\le n-1}(\P^1)\sub\Coh^{\le 2n-2}_G(X)$. Hence, the subcategories $\MM_n$ 
and $\CC^{\th}_G(X)$ have
a large intersection  provided $2n-2<\th<2n$, i.e. $n=[\th/2]+1$. 
Let us show that in this case these categories are related by tilting.

\begin{prop}\label{tilt-th-prop} 
Set $n=[\th/2]+1$. Then we have the following torsion quadruple in $\CC^{\th}_G(X)$:
\begin{equation}\label{tilt-th-eq}
(\Coh^{<\th}_G(X)[1], \TT_0, \pi^*\Coh^{\ge n}(\P^1),
[\pi^*\Coh^{\le n-1}(\P^1),\TT_1]\cap\Coh^{>\th}_G(X)).
\end{equation}
Furthermore, we have
\begin{equation}\label{tilt-th-eq2}
\MM_n=[[\pi^*\Coh^{\le n-1}(\P^1)[1],\TT_1[1]]\cap\Coh^{>\th}_G(X)[1],
\Coh^{<\th}_G(X)[1], \TT_0, \pi^*\Coh^{\ge n}(\P^1)].
\end{equation}
\end{prop}

\Pf . First, we observe that
$$(\TT_0,\pi^*\Coh^{\ge n}(\P^1), [\pi^*\Coh^{\le n-1}(\P^1),\TT_1]\cap\Coh^{>\th}_G(X))$$
is a torsion triple in $\Coh^{>\th}_G(X)$.
Indeed, this follows immediately from Proposition \ref{mod-tilt-prop} and from the
fact that the subcategory $\Coh^{>\th}_G(X)\sub\Coh_G(X)$ is stable under 
passing to quotients. Substituting this triple into the standard torsion pair
$(\Coh^{<\th}_G(X)[1],\Coh^{>\th}_G(X))$ we obtain the torsion quadruple
\eqref{tilt-th-eq}. It remains to check that all the constituents in the RHS of \eqref{tilt-th-eq2}
belong to $\MM_n$. For most of them this follows from \eqref{mod-tilt-eq}.
The remaining inclusion $\Coh^{<\th}_G(X)[1]\sub\MM_n$ is implied by the
fact that $\pi^*\OO_{\P^1}(n)$ and $\pi^*\OO_{\P^1}(n+1)$ have slope $\ge 2n>\th$.
 \ed

In conclusion we are going 
to interpret the torsion pair in $\MM_n$ arising in the above proposition in
terms of right modules over the algebra $\La(\ov{\la},\ov{m})$ (see Corollary 
\ref{quiver-cor}) assuming that $X$ is an elliptic curve. 

\begin{thm}\label{tilting-thm} Assume that $\pi:E\to\P^1$ is a ramified Galois covering with
the Galois group $G$, where $E$ is an elliptic curve, and let $\ov{m},\ov{\la}$ be
the associated ramification data.
Fix an irrational number $\th$ and set $n=[\th/2]+1$. Let us define full subcategories 
$\TT^\th,\FF^\th\sub\mod-\La(\ov{\la},\ov{m})$ as follows:
$\TT^\th$ (resp., $\FF^\th$) consists of all modules $M\simeq\oplus_{i=1}^k M_i$, where $M_i$
are indecomposable and $\deg_n(M)-\th\rk(M)<0$ (resp., $\deg_n(M)-\th\rk(M)>0$).
Then $(\TT^\th,\FF^\th)$ is a torsion pair in $\mod-\La(\ov{\la},\ov{m})$ and
one has
$$\Phi_n(\CC^{\th}_G(E))=[\FF^{\th},\TT^{\th}[-1]]\sub D^b(\La(\ov{\la},\ov{m})).$$
\end{thm}

\Pf . From Proposition \ref{tilt-th-prop} we know that $\CC^{\th}_G(E)=[\FF,\TT[-1]]$ for the
torsion pair $(\TT,\FF)$ in $\MM_n=\Phi_n^{-1}(\mod-\La(\ov{\la},\ov{m}))$ given by
$$\TT=[\pi^*\Coh^{\le n-1}(\P^1)[1],\TT_1[1]]\cap\Coh^{>\th}_G(E)[1],\
\FF=[\Coh^{<\th}_G(E)[1], \TT_0, \pi^*\Coh^{\ge n}(\P^1)].$$
We claim that one has $\Ext^1_{\MM_n}(F,T)=0$ for every $F\in\FF$ and $T\in\TT$.
It suffices to check that $\Hom_G(F,T[1])=0$ for $T\in\TT$ in the following three
cases: (i) $F\in\Coh^{<\th}_G(E)[1]$; (ii) $F\in\TT_0$; (iii) $F\in\pi^*\Coh^{\ge n}(\P^1)$.
Note that in cases (ii) and (iii) this is clear since cohomological dimension of $\Coh_G(E)$
is equal to $1$. In case (i) we obtain by Serre duality (using triviality of $\om_E$)  
$$\Hom(F,T[1])^*\simeq\Hom(T,F)=0,$$
since $T\in\Coh^{>\th}(E)[1]$ and $F\in\Coh^{<\th}(E)[1]$.

It follows that every indecomposable object of $\CC^{\th}_G(E)$ is contained either in $\TT$ or
in $\FF$. Therefore,
$\TT$ (resp., $\FF$) coincides with the full subcategory of objects $F$ such that
$F\simeq\oplus_{i=1}^n F_i$, where $F_i$ are indecomposable objects and $F_i\in\TT$
(resp., $F_i\in\FF$). Since $\deg(C)-\th \rk(C)>0$ for $C\in\CC^{\th}(E)$, it follows that
$\deg(F)-\th\rk(F)>0$ for $F\in\FF$ and $\deg(T)-\th\rk(T)<0$ for $T\in\TT$.
Taking into account Lemma \ref{deg-rk-lem} we derive that
$\TT^{\th}=\Phi_n(\TT)$ and $\FF^{\th}=\Phi_n(\FF)$.
\ed

%The set of classes of indecomposable objects maps to a certain subset $\Sigma\sub\Z^6$.
%As is observed in \cite{Sch}, $\Sigma$ can be identified with the set of roots of the affine
%Lie algebra of type $D_4^{(1)}$.

\end{document}